\documentclass[12pt]{amsart}
\usepackage{amscd,amssymb,amsthm,amsmath,enumerate,textcomp}
\usepackage[matrix,arrow]{xy}
\usepackage{url}
\usepackage{cmap}
\usepackage{hyperref}

\theoremstyle{definition}
\newtheorem{example}[equation]{Example}
\newtheorem*{example*}{Example}
\newtheorem{definition}[equation]{Definition}
\newtheorem{theorem}[equation]{Theorem}
\newtheorem{lemma}[equation]{Lemma}
\newtheorem{corollary}[equation]{Corollary}
\newtheorem{proposition}[equation]{Proposition}
\newtheorem{conjecture}[equation]{Conjecture}
\newtheorem*{conjecture*}{Conjecture}
\newtheorem*{maintheorem*}{Main Theorem}
\newtheorem*{corollary*}{Corollary}
\newtheorem{question}[equation]{Question}
\newtheorem*{question*}{Question}
\newtheorem{problem}[equation]{Problem}
\newtheorem*{problem*}{Problem}
\newtheorem*{theorem*}{Theorem}

\theoremstyle{remark}
\newtheorem{remark}[equation]{Remark}
\newtheorem*{remark*}{Remark}

\newtheorem{compactification}[equation]{Compactification Construction}

\newcommand{\PP}{{\mathbb P}}

\newcommand{\ZZ}{{\mathbb Z}}
\newcommand{\CC}{{\mathbb C}}
\newcommand{\QQ}{{\mathbb Q}}
\newcommand{\RR}{{\mathbb R}}
\newcommand{\TT}{{\mathbb T}}
\newcommand{\Aff}{{\mathbb A}}
\newcommand{\cO}{\mathcal{O}}

\makeatletter\@addtoreset{equation}{section} \makeatother

\usepackage[dvipsnames]{xcolor}
\usepackage{graphicx}

\author{Victor Przyjalkowski}

\title{On singular log Calabi--Yau compactifications of Landau--Ginzburg models
}

\thanks{This work is supported by the Russian Science Foundation under grant \textnumero 19-11-00164.
}

\address{\emph{Victor Przyjalkowski}
\newline
\textnormal{Steklov Mathematical Institute of Russian Academy of Sciences, 8 Gubkina street, Moscow, Russia.}
\newline
\textnormal{\texttt{victorprz@mi-ras.ru, victorprz@gmail.com}}}

\begin{document}

\begin{abstract}
We consider the procedure that constructs log Calabi--Yau compactifications of weak Landau--Ginzburg models of Fano varieties.
We apply it for del Pezzo surfaces and coverings of projective spaces of index one.
For the coverings of degree greater then $2$ the log Calabi--Yau compactification is singular; moreover, no
smooth projective log Calabi--Yau compactification exists.
We also prove in the cases under consideration 
the conjecture saying that the number of components of the fiber over infinity  
is equal to the dimension of an anticanonical system of the Fano variety. 
\end{abstract}

\maketitle

\section{Introduction and setup}

The mirror dual object to a smooth Fano variety $X$ is a so called \emph{Landau--Ginzburg model}, that is a quasi-projective variety $Y$ with complex-valued
regular func\-ti\-on~$\mbox{$w\colon Y\to \CC$}$ satisfying certain properties and called \emph{superpotential}.
The dual object for a given Fano variety for all versions of mirror correspondence is expected to be the same
(with different enhancements specific to the certain mirror conjectures). Homological Mirror Symmetry~\cite{Kon94} deals with singularities
of fibers of multipotentials, thus it needs the fibers to be compact (or, at least, to contain all singularities).
So it is important to have $w$ proper.
However in practice it is more convenient to have $Y$ as simple as possible and then to construct its
relative compactification. Following Givental~\cite{Gi97} and Hori--Vafa~\cite{HV00}, 
let $Y=(\CC^*)^n$, where $n=\dim (X)$. In this case the superpotential,
the regular function on the algebraic torus, can be represented by a Laurent polynomial $f$.
We associate with it a family $\left\{f=\lambda,\ \lambda\in \CC\right\}$ of fibers of the map $(\CC^*)^n\to \CC$ given by $f$.
We are interested in periods of this family.

Strong version of Mirror Symmetry conjecture of variations of Hodge structures associate \emph{a toric Landau--Ginzburg} with $X$.
The base for this conjecture is a classical Mirror Symmetry of variations of Hodge structures, see, for instance,~\cite{Gi97}.
Let $\mathbf 1$ be the fundamental class of $X$ and let $\mathcal{K}\subset H_2(X,\ZZ)$ be the set of classes of effective curves.
The series
$$
\widetilde{I}^X_{0}=
\sum_{a\in \ZZ_{\ge 0},\beta\in \mathcal{K}} \left(-K_X\cdot \beta\right)!\langle\tau_{a}
{\mathrm 1}\rangle_\beta \cdot t^{-K_X\cdot \beta}\in \CC[[t]],
$$
where $\langle\tau_{a} \mathbf 1\rangle_{\beta}$ is a \emph{one-pointed genus $0$ Gromov--Witten invariant with descendants},
see~\cite[VI-2.1]{Ma99},
is called \emph{a constant term of regularized $I$-series} 
for $X$.

Denote a constant term of a Laurent polynomial $f$ by $[f]$. Define $I_f=\sum [f^i]t^i$.

\begin{theorem}[{\cite[Proposition 2.3]{Prz08}}]
\label{Proposition: Picard--Fuchs}
Consider
a Laurent polynomial
$$
f\in \CC[x_1^{\pm 1},\ldots,x_n^{\pm 1}].
$$
Then there exist a (punctured) neighborhood $U\subset\CC$ of the point $\infty=\PP^1\setminus\CC$, a local coordinate
at $\infty$ on $U\cup\{\infty\}\subset \PP^1$, a
fiberwise $(n-1)$-form
$$
\omega_t\in \Omega^{n-1}_{(\CC^*)^n/\CC}(U),
$$
and a fiberwise $(n-1)$-cycle $\Delta_t$ over $U$
such that the Taylor expansion of the analytic function
$$
\int_{\Delta_t}\omega_t
$$
is given by the series $I_f$.
\end{theorem}

\begin{definition}
\label{definition:weak LG}
Let $X$ be a smooth $n$-dimensional Fano variety. A Laurent polynomial~$f_X$ in $n$ variables
is called a \emph{weak Landau--Ginzburg model} for $X$ if $\widetilde{I}^{X}_{0}=I_{f_X}$.
\end{definition}

\begin{remark}
\label{remark:canonical}
Mirror symmetry associates Landau--Ginzburg model not just to a Fano variety,
but to a Fano variety together with a fixed divisor class (or symplectic form) on it.
In our definition the divisor class is the anticanonical one.
\end{remark}

{
The notion of weak Landau--Ginzburg model is motivated by the following result, which can be considered as a cornerstone
of the field. Let $X$ be a smooth toric variety, let~\mbox{$\mathcal N \simeq \ZZ^N$}, $\mathcal N_\RR=\mathcal N\otimes \RR$,
let $\Sigma$ in $\mathcal N_\RR$ be its fan, and let $\{v_i\in \mathcal N|\ i\in I\}$ be
a set of integral generators of rays of $\Sigma$, where $I$ is a finite set indexing rays. We use the standard notation $x^{v}=x_1^{v^1}\cdot \ldots \cdot x_n^{v^n}$
for formal variables $x_i$ and $v=(v^1,\ldots,v^n)\in \mathcal N$.

\begin{theorem}[{\cite{Gi97}, see also~\cite{HV00}}]
\label{theorem: Givental smooth toric}
The Laurent polynomial $f_X=\sum_{i\in I}x^{v_i}$ is a weak Landau--Ginzburg model for $X$.
\end{theorem}

As we mentioned above, our goal is to get a proper family. Moreover, general mirror symmetry
expectation is that fibers of Landau--Ginzburg models are Calabi--Yau varieties mirror
dual to anticanonical sections of the initial Fano variety.
This justifies the following.

\begin{definition}
\label{definition: CY compactification}
Let $X$ be a smooth $n$-dimensional Fano variety and let $f_X$ be its weak Landau--Ginzburg model.
Consider the diagram
\begin{equation*}
\label{equation:KKP}
\xymatrix{
(\mathbb C^*)^n\ar@{^{(}->}[rr]\ar@{->}[d]_{f_X}&&Y\ar@{->}[d]^{w}\ar@{^{(}->}[rr]&&Z\ar@{->}[d]^{u}\\
\mathbb{C}\ar@{=}[rr]&&\mathbb{C}\ar@{^{(}->}[rr]&&\mathbb{P}^1,}
\end{equation*}
where $Y$ is a smooth open Calabi--Yau variety, $Z$ is smooth, $-K_Z\sim u^{-1}(pt)$, and $w$ and~$u$ are proper.
The pair $(Y,w)$ is called \emph{Calabi--Yau compactification} of $f_X$, while the pair~\mbox{$(Z,u)$} is called its \emph{log Calabi--Yau compactification}.
\end{definition}

Weak Landau--Ginzburg models for Fano threefolds are constructed, see, for instance,~\cite{Prz18}.
The way to construct log Calabi--Yau compactifications for most of them (cf. Compactification Construction~\ref{compactification construction}
below) leads to projective compactifications, see
~\cite{Prz17},~\cite{CP18}, and~\cite{ChP20}.
The projective compactifications for Fano threefolds are not smooth if and only if the threefolds lie in the family 
$X_{2.1}$ or $X_{10.1}$;
here we use the standard the notation for families of smooth Fano threefolds given in \cite{IP99}.
Note that these families are exactly ones that have base set in the anticanonical linear system, see~\cite{JR06}.
In these cases projective compactifications are $\QQ$-factorial three-dimensional varieties with ordinary double points.
Thus, they admit small non-projective resolutions. Note that the notion of tame compactified Landau--Ginzburg
model from~\cite{KKP17} used to formulate Katzarkov--Kontsevich--Pantev conjectures requires smoothness and projectivity;
however the results from~\cite{KKP17} and~\cite{Ha16} hold for the Moishezon case as well.

The phenomena to be considered in this paper is singular log Calabi--Yau compactifications.
We expect that is some cases we can not avoid singularities.

\begin{definition}
\label{definition: singular log CY}
Let $X$ be a smooth $n$-dimensional Fano variety and let $f_X$ be its weak Landau--Ginzburg model.
Consider the diagram
\begin{equation*}
\label{equation:KKP-sing}
\xymatrix{
(\mathbb C^*)^n\ar@{^{(}->}[rr]\ar@{->}[d]_{f_X}&&Y\ar@{->}[d]^{w}\ar@{^{(}->}[rr]&&Z\ar@{->}[d]^{u}\\
\mathbb{C}\ar@{=}[rr]&&\mathbb{C}\ar@{^{(}->}[rr]&&\mathbb{P}^1,}
\end{equation*}
where $Y$ is a smooth open Calabi--Yau variety, $-K_Z\sim u^{-1}(pt)$, and $w$ and~$u$ are proper.
The pair~\mbox{$(Z,u)$} is called \emph{singular log Calabi--Yau compactification} of $f_X$.
\end{definition}

In other words, Definition~\ref{definition: singular log CY} differs from Definition~\ref{definition: CY compactification}
only in allowing $Z$ to be singular over infinity.
The challenging problem is to extend the setup of Katzarkov--Kontsevich--Pantev conjectures
to singular log Calabi--Yau compactifications.

Our main interest and source of examples is weighted complete intersections.
Givental proved an analog of Theorem~\ref{theorem: Givental smooth toric} for complete intersections 
in smooth toric
varieties 
provided they admit so called \emph{nef-partitions}. (Below we give the definition applied to the case of weighted complete intersections we are interested in.)
The total space of Landau--Ginzburg model is not an algebraic torus in this case, but a complete intersection therein.
This approach was extended to the case for complete intersections in Grassmannians~\cite{BCFKS97} (see also~\cite{EHX97}) or partial flag manifolds~\cite{BCFKS98}
using ``good'' (so called small) degenerations of the ambient varieties to toric varieties and replacing the ambient variety
by the toric one.

In~\cite{Ba97} it was suggested that Theorem~\ref{theorem: Givental smooth toric} holds for a Fano variety admitting  small toric degeneration.
In a lot of cases (see the references after Definition~\ref{definition:toric condition} below) Givental's Landau--Ginzburg models for complete intersections
are birational to algebraic tori, and superpotentials for these Landau--Ginzburg models supported in the torus correspond
to toric degenerations in a way similar to the one from Theorem~\ref{theorem: Givental smooth toric}.
However the toric varieties appearing in this construction can be not ``good'': they can be very singular.
Nevertheless we claim that this approach, after some modification, can be applied for \emph{any} toric degeneration,
see, for instance,~\cite{Prz13}. 
}

Recall that for a toric variety $T$ with fan $\Sigma$ in $\mathcal{N}\simeq \ZZ^n$ the \emph{fan polytope} $F(T)$ is a convex hull in $\mathcal{N}_\RR\simeq \mathcal{N}\otimes \RR$ of integral generators of rays of $\Sigma$, and for a Laurent polyno\-mi\-al~$f\in \CC[\mathcal{N}]$, its \emph{Newton polytope} $N(f)$ is a convex hull in $\mathcal N_\RR$ of 
monomials of $f$ (considered as elements of $\mathcal N$).

\begin{definition}[{\cite[Definition 3.5]{Prz18}}]
\label{definition:toric condition}
Let $X$ be a smooth Fano variety of dimension $n$. A Laurent polynomial $f$ in $n$ variables is called \emph{toric Landau--Ginzburg
model} for $X$ if $f$ is a weak Landau--Ginzburg model for $X$, if $f$ admits a Calabi--Yau compactification,
and if there exists a degeneration $X\rightsquigarrow T$ of $X$ to a toric variety $T$ such that $N(f)=F(T)$. (The latter
condition is called \emph{the toric condition}.)
\end{definition}

{
Toric Landau--Ginzburg models exist for del Pezzo surfaces of degree at least $3$ and Fano threefolds (\cite{Prz13},~\cite{ILP13},~\cite{CCGK16},~\cite{Prz17}~\cite{KKPS19})
complete intersections (\cite{fanosearch},~\cite{ILP13},~\cite{PS15},~\cite{P16}); some partial
results are known for Grassmannians (\cite{PSh14},~\cite{PSh15}), Lagrangian Grassmannians (\cite{PR21}), and complete intersections in toric varieties
(\cite{Gi97},~\cite{DH15}). 
In~\cite[Conjecture 36]{Prz13} the existence of toric Landau--Ginzburg model for any smooth Fano variety
is declared.
}

Let us discuss what is known for smooth Fano complete intersections in weighted projective spaces
(weighted complete intersections). For a background of weighted complete intersections one can consult~\cite{Do82},~\cite{IF00},
or~\cite{PSh}, see also~\cite{PSh20-corrupt}.

Let $X$ be a smooth Fano weighted complete intersection of multidegree $(d_1,\ldots, d_k)$ in a weighted projective space
$\PP(a_0,\ldots,a_N)$.
Put $d_0=\sum a_i-\sum d_j$.

\begin{definition}
\label{definition:nef-partition}
A \emph{nef-partition} $(I_1,\ldots,I_k)$ for $X$ is a splitting
$$
\{0,\ldots,N\}=I_0\sqcup I_1\sqcup\ldots\sqcup I_k
$$
such that $\sum_{j\in I_i} a_j=d_i$ for every $i=0,\ldots,k$.
The nef-partition is called \emph{nice} if there exists an index $r\in I_0$ such that $a_r=1$.
\end{definition}

Conjecturally, nice nef-partitions exist for all smooth Fano weighted complete intersections. The evidence for this is the following.

\begin{theorem}[{\cite[Theorem 9]{Prz11}, \cite[Theorem 1.3]{PSh19}, \cite[Proposition 11.4.5]{PSh}}]
\label{theorem:nef-partitions}
Nice nef partitions exist for smooth Fano weighted complete intersections of Cartier divisors, or of codimension $2$, or of dimension at most $5$.
\end{theorem}

Given a nef-partition, in a spirit of Givental's suggestion for complete intersections in smooth toric varieties,
one can write down Landau--Ginzburg model for the weighted complete intersection.
Moreover, if the nef-partition is nice, the total space of this model is birational to an algebraic torus.

\begin{definition}
\label{definition:Givental LG}
Let $X$ be a well formed Fano weighted complete intersection.
Suppose that there exists a nice nef-partition $(I_1,\ldots,I_k)$ for $X$.
Let $a_{i0},\ldots,a_{im_i}$ denote the elements of $I_i$ for $0\leq i \leq k$, and suppose that $a_{00}=1$.
The \emph{Landau--Ginzburg model of Givental's type} is the Laurent polynomial
\begin{equation}\label{eq:fX}\tag{$\bigodot$}
f_X=\frac{(1+x_{1,1}+\ldots+x_{1,m_1})^{d_1}\cdot\ldots\cdot(1+x_{k,1}+\ldots+x_{k,m_k})^{d_k}}{\prod\limits_{i=0}^k \prod\limits_{j=1}^{m_i}
x_{i,j}^{a_{i,j}}}+x_{0,1}+\ldots+x_{0,m_0}
\end{equation}
in variables $x_{i,j}$, $1\le i \le k$, $1\le j\le m_i$.
\end{definition}


One can combinatorially check that the following holds.

\begin{proposition}
\label{proposition:weighted constant terms}
Let $X$ be a well formed Fano weighted complete intersection.
Suppose that there exists a nice nef-partition for $X$, and let
$f_X$ be defined by~\eqref{eq:fX}. Then
$$
\sum [f_X^u]t^u=\sum_j \frac{\prod_{r=0}^k(d_rv)!}{\prod_{0\leq i\leq k}^{0\leq j\leq m_i}(a_{ij}v)!}t^{d_0v}.
$$
\end{proposition}

One can deduce the following from~\cite[Theorem 0.1]{Gi97}.

\begin{theorem}[{see \cite[Theorem 1.1]{Prz07}}]
\label{theorem:weighted I-series}
Let $X$ be a smooth well formed Fano weighted complete intersection.
Suppose that there exists a nice nef-partition for $X$, and let
$f_X$ be defined by~\eqref{eq:fX}. Then
$$
\widetilde{I}_0^X=\sum_j \frac{\prod_{r=0}^k(d_rj)!}{\prod_{s=0}^N(a_sj)!}t^{d_0j}.
$$
\end{theorem}

\begin{corollary}
\label{corollary:weak LG weighted}
Let $X$ be a smooth well formed Fano weighted complete intersection.
Suppose that there exists a nice nef-partition for $X$.
Then the polynomial $f_X$ given by~\eqref{eq:fX} is a weak Landau--Ginzburg model for $X$.
\end{corollary}


{\bf
\begin{theorem}[{\cite[Theorem 3.1]{ILP13}}]
\label{theorem:weighted toric degeneration}
Let $X$ be a well formed Fano weighted complete intersection.
Suppose that there exists a nice nef-partition for $X$.
Then 
the polynomial $f_X$ given by~\eqref{eq:fX} satisfies the toric condition. 
\end{theorem}
}

The existence of Calabi--Yau compactification is proven only for complete intersections
in usual projective spaces.

{
\begin{theorem}[{\cite[Theorem 5.1]{PS15}}]
\label{theorem: CY for CI}
Let $X\subset \PP^N$ be a smooth Fano complete intersection.
Then the polynomial $f_X$ given by~\eqref{eq:fX} admits a Calabi--Yau compactification.
\end{theorem}

The proof of Theorem~\ref{theorem: CY for CI} is not so easy to generalize to weighted complete intersections.
More promising approach is the log Calabi--Yau compactification procedure suggested in~\cite{Prz17}.
First, it does not use the specificity of complete intersections. Second, it gives a precise description (up to codimension one)
of the fiber over infinity as a boundary divisor of some toric variety.
This fiber plays a key role in Katzarkov--Kontsevich--Pantev conjectures~\cite{KKP17} and in the mirror
$P=W$ conjecture~\cite{HKP20}.
This approach was applied for Fano threefolds~\cite{Prz17} and complete intersections in usual projective spaces~\cite{P16}.
Briefly, the approach is the following.
}

\begin{compactification}[{see~\cite{Prz17},~\cite{Prz18},~\cite[\S 9.3]{PSh}}]
\label{compactification construction}
Let us have a smooth Fano variety $X$ of dimension~$n$ and its weak Landau--Ginzburg model $f_X$
satisfying the toric condition. In particular, $X$ admits a degeneration to a toric variety $T$ such that $\Delta=N(f_X)=F(T)\subset \mathcal N_\RR$.
Let
$$
\nabla=\Delta^\vee=\left\{x\in \mathcal M_\RR=\mathcal N^\vee\otimes \RR \mid \langle x,y\rangle\leq -1, y\in \Delta   \right\}
$$
be the polytope dual to $\Delta$. Assume that $\nabla$ is integral (that is, $\Delta$ is reflexive) and its boundary admits a unimodular triangulation. This means that there exists a smooth toric variety $\widetilde{T}^\vee$ such that $F(\widetilde{T}^\vee)=\nabla$.
Thus $\{f_X=\lambda,\ \lambda\in \CC\}$ can be compactified to an anticanonical pencil in $\widetilde{T}^\vee$.
If $f_X$ has a particular shape, then the base set of this pencil is a union of smooth codimension $2$ subvarieties.
Moreover, the base set does not contain codimension $2$ toric strata of $\widetilde{T}^\vee$.
If we blow up one of the components of the base set, we get a 
pencil whose base set is a union of smooth codimension $2$ subvarieties again,
and the anticanonical class of the blown up family is still linearly equivalent to a fiber. Repeating such blow ups one by one, we arrive to the required log Calabi--Yau compactification $u\colon Z\to \PP^1$.
As a bonus of this procedure we get a description of the fiber over infinity $u^{-1}(\infty)$
as the boundary divisor (up to codimension~$2$) for $\widetilde{T}^\vee$. In particular, it consists of $k$ components, where $k$ is the number
of integral boundary points of $\nabla$.
\end{compactification}

This log Calabi--Yau compactification procedure suggests the following.

\begin{conjecture}[{cf. \cite[Conjecture 1.6]{ChP20}}]
\label{conjecture:linear systems}
Let $X$ be a smooth Fano variety, and let~$(Z,u)$ be a singular log Calabi--Yau compactification of its toric Landau--Ginzburg model.
Then $u^{-1}(\infty)$ consists of
$$
\chi\big(\mathcal{O}(-K_X)\big)-1=h^0\big(\mathcal{O}_X(-K_{X})\big)-1
$$
irreducible components.
\end{conjecture}

Conjecture~\ref{conjecture:linear systems} is proven by direct calculations for Fano threefolds,
complete intersections in projective spaces, and in some toric varieties, see~\cite{ChP20}.
Surprisingly, it is not checked for del Pezzo surfaces; we fill this gap in Proposition~\ref{proposition:del Pezzo}.

\begin{remark}
\label{remark:linear system}
Conjecture~\ref{conjecture:linear systems} is motivated by the following.
Let $\Delta$ be a Newton polytope for $f_X$. Consider
the degeneration of $X$ to the toric variety $T_\Delta$. Since this degeneration is flat, one has
$$
\chi(X)=\chi(T_\Delta).
$$
On the other hand, the log pair $(T_\Delta,-K_{T_\Delta})$ is Kawamata log terminal, see \cite[Proposition~3.7]{Ko95} for the terminology.
Thus, by Kodaira vanishing (see, for example,~\cite[Theorem~2.70]{KM98}),
one has
$$
h^0(-K_X)=h^0(-K_{T_\Delta}).
$$
The anticanonical linear system of $T_\Delta$ can be described as a linear system
of Laurent polynomials supported on
the dual polytope $\nabla$.
Suppose that the log Calabi--Yau compactification procedure described in Compactification Construction~\ref{compactification construction}
is applicable for~$f_X$.
In particular, $\nabla$ is integral and $T_\nabla$ admits a toric crepant resolution
$\widetilde T_\nabla\to T_\nabla$.
The dimension of the anticanonical linear system of $T_\Delta$ is the number of integral points on the boundary of $\nabla$.
Since these boundary points are in one-to-one correspondence with boundary divisors of $\widetilde T_\nabla$
and, thus, with irreducible components of the fiber~\mbox{$u^{-1}(\infty)$}, the conjecture follows.
\end{remark}

The fiber over infinity is not the only one that can be reducible. Another reducible fibers of log Calabi--Yau compactification
also predict invariants of the initial Fano variety.

\begin{conjecture}[{\cite[Conjecture 1.1]{PS15}}]
\label{conjecture:mirror Hodge}
Let $X$ be a smooth Fano variety of dimension $n$ and let $(Y,w)$ be a Calabi--Yau compactification of a toric Landau--Ginzburg model for $X$.
Put
\begin{multline*}
\kappa_{Y}= \sharp\mbox{\big[irreducible components of all reducible
fibers of $(Y,w)$\big]}
- \sharp \mbox{\big[reducible fibers\big]}.
\end{multline*}
Then one has $h_{pr}^{1,n-1}(X)=\kappa_{Y}$, where $h_{pr}^{1,1}(X)=h^{1,1}(X)+1$ and $h_{pr}^{1,n-1}(X)=h^{1,n-1}(X)$ for $n>2$.
\end{conjecture}

Compactification Construction~\ref{compactification construction} enables one to verify this conjecture, see Example~\ref{example:double sextic}. However usually
it is more convenient to use more straightforward compactification in a projective space (see~\cite{CP18}) for this.

Unfortunately, Compactification Construction~\ref{compactification construction} has some restrictions.
First, we need $\Delta$ to be reflexive. Then, we assume that $\nabla$ admits an unimodular triangulation, so that there
is a smooth toric variety whose fan polytope is $\nabla$. Finally, we need to have ``good'' base set of the pencil compactified in
a smooth toric variety.

In fact toric Landau--Ginzburg models usually have very specific Newton polytopes and coefficients,
so that the second and third problems are relatively easy to solve. However one can not launch the compactification procedure
as it is for the non-reflexive case, which is common already for weighted complete intersections.

We claim that Compactification Construction~\ref{compactification construction} after some modification can be applied
in more general cases. This is done in the following main results of the paper.

\begin{proposition}[Proposition~\ref{proposition:del Pezzo}]
Let $S$ be a smooth del Pezzo surface. Then standard (see Section~\ref{proposition:del Pezzo}) weak Landau--Ginzburg model $f_S$ for $S$ admits a log Calabi--Yau compactification.
In particular, $f_S$ is a toric Landau--Ginzburg model. Moreover, Conjectures~\ref{conjecture:linear systems} and~\ref{conjecture:mirror Hodge} hold for it.
\end{proposition}

\begin{theorem}[Theorem~\ref{theorem:coverings}]
Let $X$ be a smooth Fano index one covering of a projective space.
Then weak Landau--Ginzburg model of Givental's type
is a toric Landau--Ginzburg model for $X$. It admits a singular log Calabi--Yau compactification.
Conjecture~\ref{conjecture:linear systems} holds for it.
\end{theorem}

In Lemma~\ref{lemma: non-projective} we show that there is no smooth projective log Calabi--Yau compactification in the case
considered in Theorem~\ref{theorem:coverings} if the degree of the covering is greater then $2$;
we expect that there is no smooth non-projective compactification as well.


In Section~\ref{section:discussion} we outline the suggested generalization of the compactification procedure.

\smallskip
The author is grateful to
I.\,Cheltsov, A.\,Corti, T.\,Pantev, and C.\,Shramov for useful discussions.

\section{Del Pezzo surfaces}
\label{section:del Pezzo}

Recall that smooth del Pezzo surfaces can be described as either a quadric surface or an anticanonical degree $d\leq 9$ surface $S_d$ which is a blow up of $\PP^2$ in $9-d$
general enough points. For $d\geq 3$ the anticanonical class of $S_d$ is very ample.
Landau--Ginzburg models for smooth del Pezzo surfaces and general divisors on them were suggested in~\cite{AKO06}.
They are pencils of elliptic curves with nodal singular fibers and one fiber (``over infinity'') which is a wheel of $12-d$
smooth rational curves. Construction of toric Landau--Ginzburg model for \emph{any} divisor on $S_d$, $d\geq 3$, one can find in~\cite{Prz17}.
This construction uses Gorenstein toric degenerations obtained in~\cite{Ba85}.
For general divisors log Calabi--Yau compactifications of the toric Landau--Ginzburg models are pencils constructed in~\cite{AKO06}.
For simplicity we consider toric Landau--Ginzburg corresponding to the anticanonical divisor;
note that in fact the only thing in the proof of Proposition~\ref{proposition:del Pezzo} that depend on the particular divisor
is the structure of the central fiber.
Del Pezzo surfaces of degrees $1$ and $2$ are weighted complete intersections, so they have weak Landau--Ginzburg
models of Givental's type. We call the weak Landau--Ginzburg models mentioned above (for all del Pezzo surfaces)
\emph{standard}.

\begin{proposition}
\label{proposition:del Pezzo}
Let $S$ be a smooth del Pezzo surface. Then standard weak Landau--Ginzburg model $f_S$ for $S$ admits a log Calabi--Yau compactification.
In particular, $f_S$ is a toric Landau--Ginzburg model. Moreover, Conjectures~\ref{conjecture:linear systems} and~\ref{conjecture:mirror Hodge} hold for it.
\end{proposition}

\begin{proof}
Note that toric condition holds for standard weak Landau--Ginzburg models of del Pezzo surfaces by definition.
Thus to prove the second assertion of the proposition it is enough to prove the first one.
Let $\deg (S)\ge 3$. Then $f_{S}$ admits the log Calabi--Yau compactification given by~\cite[Compactification Construction 16]{Prz17}.
Conjecture~\ref{conjecture:linear systems} holds for $S$ by~\cite[Remark 17]{Prz17}.
One can check Conjecture~\ref{conjecture:mirror Hodge} compactifying $f_{S}$ as a pencil of elliptic curves in $\PP^2$ and resolve
the base set of the pencil. (In fact this is just another description of Compactification Construction~\ref{compactification construction}.)
For instance, if $d=3$, then
$$
f_S=\frac{(x+y+1)^3}{xy}.
$$
Using the embedding $\TT(x,y)\subset \PP(x:y:z)$,
where~\mbox{$\TT(x,y)\cong(\CC^*)^2$} denotes the algebraic torus
with coordinates $x$ and $y$,
one can compactify the family of fibers for~$f_S$ to the elliptic pencil
$$
\mathcal F=\{\mathcal F_{(\lambda:\mu)}\mid (\lambda:\mu)\in \PP^1\},
$$
where
$$
\mathcal F_{(\lambda:\mu)}=
\{\mu (x+y+z)^3=\lambda xyz\}.
$$
This pencil is generated by the triple line
$$
l=\{x+y+z=0\}
$$
and the union of lines
$$
l_1=\{x=0\},\quad l_2=\{y=0\},\quad l_3=\{z=0\}.
$$
The base set of $\mathcal F$ consists of three points $p_i=l\cap l_i$, $i=1,2,3$. At each of them, the fiber over the point $(0:1)\in\PP^1$
has multiplicity $3$, and any other fiber
has multiplicity $1$.
Let us resolve the base set of $\mathcal F$. First, resolve
it in a neighborhood of the point $p_1$.
For this blow it up. Since the multiplicity of $\mathcal{F}_{(0:1)}=3l$ at $p_1$ is $3$,
the exceptional curve lies in the fiber over $(0:1)$ of the proper transform ${\mathcal{F}}^1$ of the pencil $\mathcal{F}$,
and the multiplicity of this fiber along the curve is $2$. So the intersection of the exceptional
curve with the fiber over infinity is a base point of~${\mathcal{F}}^1$ of multiplicity $2$. Proceeding like this, make two more blow ups to resolve the base set in the neighborhood of $p_1$; note that the first two exceptional divisors lie in the fiber over~\mbox{$(0:1)$} of the resolved pencil,
while the last one is horizontal with respect to it, that is, it projects surjectively to the base of the pencil.
Similarly, resolve the base set of $\mathcal F$ in a neighborhood of the points $p_2$ and $p_3$.

Note that after this resolution the anticanonical class is linearly equivalent to a fiber of the proper transform of~$\mathcal{F}$,
so we arrive to a log Calabi--Yau compactification of $f_S$.
Its fiber over $(1:0)$ consists of three curves, namely, of the proper transforms of $l_1$, $l_2$, and~$l_3$.
Therefore, Conjecture~\ref{conjecture:linear systems} holds for~$S$.

Finally, we point out that there are three exceptional curves over each of the three points~$p_i$.
Two of them lie in the fiber over the point $(0:1)$
(i.e., in the fiber containing the triple line),
and one is horizontal. Hence we get
$$
1+3\cdot 2=7
$$
irreducible components in the fiber over $(0:1)$, so~\ref{conjecture:mirror Hodge} holds in this case. These components form the extended
Dynkin diagram~$\widetilde{\mathbf{E}}_6$.

Now let $\deg (S)=2$. Then $S$ can be described as a quartic hypersurface in $\PP(1,1,1,2)$. Thus, Givental's Landau--Ginzburg model
for $S$ is the polynomial
$$
f_{2}=\frac{(x+y+1)^4}{xy}.
$$
Its Newton polytope $\Delta$ is the convex hull of vertices $(3,-1)$ , $(-1,3)$, $(-1,-1)$.
The dual polytope $\nabla=\Delta^\vee$ is the convex hull of vertices $(1,0)$, $(0,1)$, $(-\frac{1}{2},-\frac{1}{2})$.
Thus $\nabla$ is not integral, so that we can not apply Compactification Construction~\ref{compactification construction}.
However we still can construct a toric variety whose rays are generated by the vertices of $\nabla$.
Its fan is generated by integral points $(1,0)$, $(0,1)$, $(-1,-1)$, so
it is nothing but $\PP^2$. The pencil of fibers for $f_2$ is now generated by the fiber over infinity,
which is a union of three lines $l_1$, $l_2$, $l_3$, one of which has multiplicity $2$, and the central fiber,
which is a line $l$ in general position taken with multiplicity $4$.
In other words, as above we can compactify the pencil given by $f_2$ to the pencil
$\mathcal F=\{\mathcal F_{(\lambda:\mu)},\ (\lambda:\mu)\in \PP^1\}$, where
$$
\mathcal F_{(\lambda:\mu)}=\{\mu (x+y+z)^4=\lambda xyz^2\},
$$
so that
$$l=\{x+y+z=0\},\qquad l_1=\{x=0\},\qquad l_2=\{y=0\},\qquad l_3=\{z=0\}.
$$
In par\-ti\-cu\-lar,
$$\mathcal F_0=\mathcal F_{(0:1)}=4l\qquad \mbox{and}\qquad \mathcal F_\infty=\mathcal F_{(1:0)}=l_1+l_2+2l_3.
$$
We also have $\mathcal F_t\sim -K_{\PP^2}+l_3$.

Denote $p_i=l\cap l_i$. Let us resolve the base set of $\mathcal F$. First, resolve
it in a neighborhood of the point $p_1$.
For this blow it up. Since the multiplicity of $\mathcal{F}_{(0:1)}=4l$ at $p_1$
is $4$, we see that
the exceptional curve of the blow up lies in the fiber over $(0:1)$ of the proper transform ${\mathcal{F}}^1$ of the pencil $\mathcal{F}$,
and this curve has multiplicity $3$ in this fiber. So the intersection of the exceptional
curve with the fiber over $(1:0)$ is a base point of multiplicity $3$ of~${\mathcal{F}}^1$. Proceeding like this, make three more blow ups to resolve the base set of $\mathcal F$
in the neighborhood of $p_1$; note that the first three exceptional divisors lie in the fiber over $(0:1)$ of the resolved pencil,
while the last one is horizontal.
Similarly, resolve the base set of $\mathcal F$ in a neighborhood of the point $p_2$.

Now blow up the point $p_3$. Let $e$ be the exceptional curve of the latter
blow up.
Since the multiplicity of $\mathcal F_{(0:1)}$ at $p_3$ is equal to $4$, while the multiplicity of $\mathcal F_{(1:0)}$ at $p_3$ is equal to $2$, we conclude that
the multiplicity of the fiber over $(0:1)$ of the proper transform of $\mathcal{F}$ at $e$ is $4-2=2$.
Denote the proper transforms of the curves by the same symbols as the initial curves.
To resolve the base set of the pencil, we need to blow up the point~$e_1\cap l_3$.
After this blow up we arrive to a base point free family
$$
\widetilde{\mathcal F}=\{\widetilde{\mathcal F}_t\mid t\in \PP^1\}
$$
with total space $\widetilde Z$.
One has
$$
\widetilde{\mathcal F}_t\sim-K_{\widetilde{Z}}+l_3.
$$
Moreover, since $l_3^2=1$ on $\PP^2$ and we blow up
two smooth points lying on $l_2$, one has~\mbox{$l_3^2=-1$} on $\widetilde Z$.
Thus, contracting $l_3$ on $\widetilde{Z}$, we get the pencil
$\mathcal G=\{\mathcal G_t\mid t\in \PP^1\}$
with smooth total space $Z$ and $-K_Z\sim \mathcal G_t$. This means that $Z$ is the required log Calabi--Yau compactification for~$f_2$.
Conjectures~\ref{conjecture:linear systems} and~\ref{conjecture:mirror Hodge} easily follow from this construction.
Indeed, $\mathcal G_{(1:0)}=l_1+l_2$ is a union of two rational curves intersecting transversally at two points. Note that the fiber $\mathcal G_{(0:1)}$ consists of
$l$, $e$, three exceptional divisors lying over $p_1$, and three exceptional divisors lying over $p_2$.
These curves form the extended
Dynkin diagram~$\widetilde{\mathbf{E}}_7$.

Similarly one can get the log Calabi--Yau compactification for $\deg (S)=1$. In this case~$S$ can be described as a sextic hypersurface in $\PP(1,1,2,3)$. Thus, Givental's Landau--Ginzburg model
for $S$ is the polynomial
$$
f_{1}=\frac{(x+y+1)^6}{xy^2}.
$$
The polytope dual to $N(f_1)$ is a convex hull of points $(1,0)$, $(0,\frac{1}{2})$, $(-\frac{1}{3},-\frac{1}{3})$.
We compactify our family to the pencil $\mathcal F=\{\mathcal F_{(\lambda:\mu)},\ (\lambda:\mu)\in \PP^1\}$, where
$$
\mathcal F_{(\lambda:\mu)}=\{\mu (x+y+z)^6=\lambda xy^2z^3\},
$$
so that for
$$l=\{x+y+z=0\},\qquad l_1=\{x=0\},\qquad l_2=\{y=0\},\qquad l_3=\{z=0\}
$$
we have
$$\mathcal F_0=\mathcal F_{(0:1)}=6l\qquad \mbox{and}\qquad \mathcal F_\infty=\mathcal F_{(1:0)}=l_1+2l_2+3l_3.
$$
We also have $\mathcal F_t\sim -K_{\PP^2}+l_2+2l_3$.

Denote $p_i=l\cap p_i$. As above, making six blow ups for $p_1$, three blow ups for $p_2$, and two blow ups for $p_3$, we get
the family $\widetilde{\mathcal F}=\{\widetilde{\mathcal F}_t,\ t\in \PP^1\}$ with total space $\widetilde{Z}$ and $\widetilde{\mathcal F}_t\sim-K_{\widetilde{Z}}+l_2+2l_3$.
Note that~$l_2^2=-2$ and~$l_3^2=-1$ on $\widetilde{Z}$. Let $\varphi\colon \widetilde Z\to Z'$ be the contraction of $l_3$.
Then $Z'$ is smooth and~$\varphi(l_2)^2=-1$, so that we can contract $\varphi(l_2)$ and get the log Calabi--Yau compactification $Z$ of $f_X$.
Conjectures~\ref{conjecture:linear systems} and~\ref{conjecture:mirror Hodge} easily follow from this construction.
Note that the fiber of the pencil $Z\to\PP^1$ over $(0:1)$ consists of~$l$,
five exceptional divisors lying over $p_1$, two exceptional divisors lying over $p_2$, and one exceptional divisor lying over~$p_3$.
These curves form the extended Dynkin diagram~$\widetilde{\mathbf{E}}_8$.
\end{proof}

\section{Index one coverings of projective spaces}
\label{section:coverings}
In this section we generalize the log Calabi--Yau compactification procedure for del Pezzo surface of degree $2$ constructed in the
proof of Proposition~\ref{proposition:del Pezzo} to Fano index one smooth Fano coverings of projective spaces.
Let $X$ be Fano index $1$ smooth Fano variety which is a $a$-to-$1$-covering of a projective space branched in a divisor of degree $ad$.
Let
$$
\alpha=(a-1)d=\dim (X).
$$
Then $X$ can be described as a hypersurface of degree $ad$ in
$$
\PP(\underbrace{1,\ldots,1}_{\alpha+1},d).
$$
Thus its Givental's Landau--Ginzburg model is
$$
f_X=\frac{(x_1+\ldots+x_{\alpha}+1)^{ad}}{x_1\cdot\ldots\cdot x_{\alpha}}.
$$

\begin{theorem}
\label{theorem:coverings}
The Laurent polynomial $f_X$
is a toric Landau--Ginzburg model for $X$. It admits a singular log Calabi--Yau compactification.
Conjecture~\ref{conjecture:linear systems} holds for it.
\end{theorem}

\begin{proof}
The polynomial $f_X$ is a weak Landau--Ginzburg model by Corollary~\ref{corollary:weak LG weighted}. It satisfies the toric condition by Theorem~\ref{theorem:weighted toric degeneration}. Now construct a Calabi--Yau compactification.

The Newton polytope $\Delta$ for $f_X$ is a convex hull of the points whose coordinates are given by the
rows of the matrix
$$
\left(%
\begin{array}{rrrr}
  ad-1 &  -1 & \ldots & -1 \\
   -1 & ad-1 & \ldots & -1 \\
  \ldots & \ldots & \ldots & \ldots \\
  -1 & -1 & \ldots & ad-1 \\
  -1 & -1 & \ldots & -1 \\
\end{array}%
\right).
$$

The dual polytope $\Delta^\vee$ is a convex hull of the points whose coordinates are given by the
rows of the matrix
$$
\left(%
\begin{array}{rrrr}
  1 & 0 & \ldots & 0 \\
  0 & 1 & \ldots & 0 \\
  \ldots & \ldots & \ldots & \ldots\\
  0 & 0 & \ldots & 1\\
  -\frac{1}{d} & -\frac{1}{d} & \ldots & -\frac{1}{d}\\
\end{array}%
\right).
$$
The toric variety whose rays are generated by rows of this matrix is $\PP=\PP^{\alpha}$, and integral generators of these rays
are the rows of the latter matrix with the last one scaled by $d$. This means that the compactified (in $\PP$) pencil $\mathcal F$ corresponding
to $f_X$ is generated by the fiber over infinity $\mathcal F_\infty$, which is a union of coordinate hyperplanes in $\PP$ with one of them taken with multiplicity $d$, and
the central (that is, lying over $(0:1)$) fiber $\mathcal F_0$, which is a hyperplane (in general position with respect to the coordinate hyperplanes) taken with multiplicity $ad$.
In other words, we compactify the family via the embedding
$$
\TT(x_1,\ldots,x_{\alpha})\subset \PP(x_0:\ldots:x_{\alpha}),
$$
where~\mbox{$\TT(x_1,\ldots,x_{\alpha})\cong(\CC^*)^\alpha$} denotes the algebraic torus
with coordinates $x_1,\ldots,x_\alpha$, to get the pencil
$\mathcal F=\{\mathcal F_{(\lambda:\mu)},\ (\lambda:\mu)\in \PP^1\}$, where
$$
\mathcal F_{(\lambda:\mu)}=\{\mu \cdot (x_0+\ldots+x_{\alpha})^{ad}=\lambda \cdot x_0^d\cdot x_1\cdot\ldots\cdot x_{\alpha}\}.
$$
In particular, for $H=\{x_0+\ldots+x_{\alpha}=0\}$ and $H_i=\{x_i=0\}$, $i=1,\ldots,\alpha$, we have
$$
\mathcal F_0=\mathcal F_{(0:1)}=adH,\ \ \ \mathcal F_\infty=\mathcal F_{(1:0)}=dH_0+H_1+\ldots+H_{\alpha}.
$$
Since $-K_\PP\sim H_0+\ldots+H_\alpha$, we have $\mathcal F_t\sim -K_{\PP}+(d-1)H_0$.
If we denote $B_i=H\cap H_i$, then the base set of $\mathcal F$ is the union of $B_i$, $i=0,\ldots,\alpha$.
The multiplicity of $\mathcal F_0$ at ${B_i}$ is equal to $ad$, while
the multiplicity of $\mathcal F_\infty$ at ${B_0}$ is equal to $d$ and the multiplicity of $\mathcal F_\infty$ at $B_i$, $i=1,\ldots,\alpha$,
is equal to $1$.

Blow up $\PP$ along $B_0$. We get an exceptional divisor $E_1$ and the proper transform  $\mathcal F^1$
of the pencil $\mathcal{F}$. One has
$\mathcal F^1_0=(a-1)dE_1+adH$ and
$$
\mathcal F^1_\infty=dH_0+H_1+\ldots+H_{\alpha}.
$$
As before, we denote divisors and their proper transforms by the same symbols.
Furthermore, denote~\mbox{$B_0=E_1\cap H_0$} and
$
B_i=H\cap H_i 
$
for $i=1,\ldots,\alpha$. The multiplicity of $F^1_0=(a-1)d$ in $B_0$ is equal to $(a-1)d$, the multiplicity of $\mathcal F^1_\infty$
at $B_0$ is equal to $d$, the multiplicity of $\mathcal F^1_0$ at $B_i$, $i=1,\ldots,\alpha$, is equal to $ad$,
and the multiplicity of $\mathcal F^1_\infty$ at $B_i$, $i=1,\ldots,\alpha$, is equal to $1$.
The base set of $\mathcal F^1$ is a union of several smooth codimension $2$ strata.
Its irreducible components are
$B_i$, $i=0,\ldots,\alpha$, and $B^1_i=E_1\cap H_i$, $i=1,\ldots,\alpha$.
We also have
$$
\mathcal F^1_t\sim-K_{Z^1}+(d-1)H_0,
$$
where $Z^1$ is the total space of $\mathcal F^1$.

Repeating this procedure $a-1$ more times, one can resolve the component $B_0$ of the base set
and get the pencil $\mathcal F^a=\{\mathcal F^a_t, \ t\in \PP^1\}$ with (smooth) total space $Z^a$. The base set of $\mathcal F^a$ is a union of smooth codimension $2$ subvarieties of $Z^a$ that intersect each other
transversally. The divisor
$\mathcal F^a_\infty$ has multiplicity $1$ at every component of the base set. Moreover, one has
$\mathcal F^a_t\sim-K_{Z^a}+(d-1)H_0$.

Choose a toric structure on the original $\PP^\alpha$ with boundary divisors
$H_0, H, H_2,\ldots, H_\alpha$. Then the blow ups we consider above are toric, and thus $Z^a$ is also toric.
The rays of the fan of $Z^a$ are generated by the vectors
$$
\begin{array}{l}
  v_0=(1,0,0,\ldots,0,0), \\
  v_1=(0,1,0,\ldots,0,0), \\
  \ldots \\
  v_{\alpha-1}= (0,0,0,\ldots,0,1),\\
  v_{\alpha}=(-1,-1,-1,\ldots,-1,-1), \\
  u_1=(1,1,0,\ldots,0,0),\\
  \ldots \\
  u_a=(a,1,0,\ldots,0,0).\\
\end{array}
$$
Here the vector $v_1$ corresponds to the proper transform of $H$,
$v_i$, $i=0,2,\ldots,v_\alpha$, correspond to proper transforms of $H_i$,
and $u_i$, $i=1,\ldots,u_a$, correspond to exceptional divisors $E_i$. 
Since 
$$
(a-1) v_0=u_a+v_2+\ldots+v_{\alpha},
$$
one can check that we can blow down the proper transform of $H_0$ on $Z^a$ to a point $P$.
(In other words, the normal bundle to $\PP^{\alpha-1}\cong H_0\subset Z^a$ is $\cO_{\PP^{\alpha-1}}(1-a)$, so we can blow strict
transform of $H_0$ down  to a point.)
The point $P$ is a cyclic quotient singularity of type $\frac{1}{a-1}(1,\ldots,1)$. In particular, it is smooth for $a=2$ and singular for $a>2$.
Let $Z^a\to Z'$ be the blow down. The pencil $\mathcal F^a$ descends to a pencil $\mathcal F'=\{\mathcal F'_t\mid t\in \PP^1\}$ on~$Z'$.
Note that the point $P$ does not lie in the base set of $\mathcal F'$, because $H_0\subset Z^a$ does not intersect $\mathcal F^a_t$ for $t\neq \infty$.
Moreover, one has $\mathcal F'_t\sim -K_{Z'}$.

Now we can proceed as in Compactification Construction~\ref{compactification construction}, cf.~\cite[Proposition 26]{Prz17},
and get a 
family $u\colon Z\to\PP^1$. The variety $Z$ is smooth for $a=2$ and has a unique singular point for $a>2$, and this point lies in $D=u^{-1}(\infty)$.
Moreover, one has $-K_{Z}\sim D$, so that $Z$ is a required singular log Calabi--Yau compactification.
In particular, $Z\setminus D$ is a Calabi--Yau compactification of $f_X$, so the latter Laurent polynomial
is a toric Landau--Ginzburg model for $X$.

The fiber over $(1:0)$ of the family $Z\to\PP^1$ 
is a normal crossing divisor, and its
dual intersection complex
is a boundary of $a(d-1)$-simplex with one contracted face of maximal dimension,
so it is homeomorphic to an $\left(a(d-1)-1\right)$-sphere and consists of $\alpha=(a-1)d$ components.
By~\cite[Theorem 3.3.4]{Do82} and~\cite[Corollary 3.3]{PrzyalkowskiShramov-AutWCI} one has
$$h^0(\cO_X(-K_X))=\alpha+1,$$
so Conjecture~\ref{conjecture:linear systems} follows.
\end{proof}

\begin{lemma}
\label{lemma: non-projective}
Let $a>2$. Then the toric Landau--Ginzburg model $f_X$ does not admit a smooth projective log Calabi--Yau compactification.
\end{lemma}

\begin{proof}
We use the notation of Theorem~\ref{theorem:coverings}.
Let $w\colon W\to \PP^1$ be a smooth projective log Calabi--Yau compactification for $f_X$.
Note that $Z'$ is $\QQ$-factorial by construction.
By~\cite[Theorem 1]{Ka08} the birational isomorphism of pairs $\left(W,w^{-1}(\infty)\right)$
and $(Z',\mathcal F_\infty')$ 
is composed from a sequence of flops.
Since $a>2$, the point $P$ is singular, so the isomorphism is not biregular.
In particular, there should exist a flop
in the neighborhood of $P$.
Using the torus action one can see that 
the flopping locus should be torus invariant.
The point $P$ correspond to a cone with generators $u_a,v_2,\ldots,v_\alpha$.
The affine hyperplane
$$
L=\left\{\frac{\alpha}{a-1}e_1+\frac{a-\alpha a-1}{a-1}e_2+ e_3+\ldots+e_\alpha=1\right\},
$$
where $e_1,\ldots,e_\alpha$ are standard basis vectors in $\ZZ^\alpha$,
contains $v_2,\ldots,v_\alpha, u_a$ and no other generators of rays of the fan of $Z'$.
This shows that there is no torus invariant flopping locus in the neighborhood of $P$.
The contradiction proves the lemma.
\end{proof}

\begin{question}
\label{question:singular log CY}
Is it true that the Laurent polynomial $f_X$ from Theorem~\ref{theorem:coverings} does not admit a smooth algebraic (non-projective)
log Calabi--Yau compactification? Can definitions and constructions from~\cite{KKP17} be generalized to this case?
\end{question}

One of approaches to positive answer on the first part of Question~\ref{question:singular log CY} is to prove
the analog of~\cite[Theorem 1]{Ka08} for non-projective case. Unfortunately, the proof of this theorem relies
on the Minimal Model Program, which is hard to be generalized to the non-projective case.

\begin{remark}[{cf.~\cite[Remark 19]{Prz13}}]
The initial Givental's Landau--Ginzburg model for $X$ is the variety
$$
Y'=\{y_1+\ldots+y_{\alpha+1}=1\}\subset (\CC^*)^{\alpha+1}
$$
with the function $w'=y_1+\ldots+y_{\alpha+1}+\frac{1}{y_1\cdot\ldots\cdot y_\alpha\cdot y_{\alpha+1}^d}$
or, after shifting by one, with the function $w=\frac{1}{y_1\cdot\ldots\cdot y_\alpha\cdot y_{\alpha+1}^d}$.
The family $Z\to \PP^1$ is a compactification of $w\colon Y'\to \CC$.
After inverting the superpotential, $t=1/w$, we get the family $Y'\to \CC$ given by~$\mbox{$t=y_1\cdot\ldots\cdot y_\alpha\cdot y_{\alpha+1}^d$}$. This family is nothing but a Landau--Ginzburg model
for
$$
\PP(\underbrace{1,\ldots,1}_{\alpha},d)
$$
in the sense of~\cite{CG11}. Thus, Lemma~\ref{lemma: non-projective} shows that this family can not be compactified
to a projective family of Calabi--Yau varieties with smooth total space.
\end{remark}

Theorem~\ref{theorem:coverings} implies the following.

\begin{corollary}
\label{corollary:Hodge-Tate}
Let $W$ be any projective log Calabi--Yau compactification of $f_X$. Then its cohomology
has a pure Hodge structure, and it is of Hodge--Tate type, that is, $h^{p,q}(W)=0$ if $p\neq q$.
\end{corollary}

\begin{proof}
Note that $Z'$ is toric, and, thus its cohomology has a pure Hodge structure of Hodge--Tate type. Since $Z$ is obtained from $Z'$ by blow ups in smooth centers of Hodge--Tate type, $Z$ is of Hodge--Tate type as well. Finally, since, as above, $W$ is obtained from $Z$ by a sequence of flops,
and Hodge--Tate property does not change under flops, the corollary follows.
\end{proof}

The construction used in the proof of Theorem~\ref{theorem:coverings} enables one to count components of the cental fiber to prove Conjecture~\ref{conjecture:mirror Hodge}.

\begin{example}
\label{example:double sextic}
Let $X\subset \PP(1,1,1,1,3)$ be a three-dimensional sextic double solid.
Let us construct the log Calabi--Yau compactification following the proof of Theorem~\ref{theorem:coverings}.
We also keep track
of exceptional divisors lying in the central fiber (which turns out the only reducible one, except for the fiber over
infinity).

One has
$$
f_X=\frac{(x_1+x_2+x_3+1)^6}{x_1x_2x_3},
$$
so that $\mathcal F=\{\mathcal F_{(\lambda:\mu)},\ (\lambda:\mu)\in \PP^1\}$, where
$$
\mathcal F_{(\lambda:\mu)}=\{\mu \cdot (x_0+x_1+x_2+x_3)^{6}=\lambda\cdot x_0^3x_1x_2x_3.
$$
Denoting
$$H=\{x_0+x_1+x_2+x_3=0\}\qquad \mbox{and}\qquad H_i=\{x_i=0\}\ \  \mbox{for}\  i=1,2,3,
$$
we have
$$\mathcal F_0=6H\qquad \mbox{and}\qquad \mathcal F_\infty=3H_0+H_1+H_2+H_3.
$$
Denote also $B_i=H\cap H_i$ 
and $P_{ij}=H\cap H_i\cap H_j$.

To resolve the base set of $\mathcal F$,
we blow up $B_0$ and the intersection of the exceptional divisor with the proper transform of $H_0$.
Then we blow down the proper transform $H_0$ to a smooth point outside the base set of the
corresponding pencil, and resolve the base set blowing up smooth curves.
Therefore, to count components of the central fiber it is enough to count exceptional divisors of
the blow ups lying over general points of $B_i$ and over $P_{ij}$.

Let us make a general useful observation. Let $\mathcal J=\{\mathcal J_t\mid t\in \PP^1\}$ be a pencil on a smooth threefold $Y$.
Let $C\subset Y$ be a smooth curve such that the multiplicity of $\mathcal J_0$ in $C$ is equal to $bm$ and the multiplicity of $\mathcal J_\infty$ in $C$ is equal to $m$ for positive integers $b$ and $m$.
Then to resolve the base locus of $\mathcal J$ in the neighborhood of a general point of $C$ we need to make $b$ blow ups of curves lying over $C$.
Moreover, one exceptional divisor for this resolution is not contained in any fiber of the proper transform of $\mathcal J$, while
$b-1$ ones lie in the fiber over $0$.
In our situation, since the multiplicity of
$\mathcal F_0$ in ${B_i}$, $i=0,\ldots,3$, is equal to $6$, the multiplicity of $\mathcal F_\infty$ in ${B_0}$ is equal to $3$,
while the multiplicity of
$\mathcal F_\infty$ in ${B_i}$, $i=1,2,3$, is equal to
$1$,
we have
$$
\left(\frac{6}{3}-1\right)+3\cdot\left(\frac{6}{1}-1\right)=16
$$
exceptional divisors in the central fiber
lying over general points of the irreducible components of the base set.

Now we need to count exceptional divisors lying over
the points $P_{ij}$.
Set~\mbox{$P=P_{0j}$} for some $1\le j\le 3$.
Let~$E$ be the exceptional divisor of the blow up of $B_0$. Then, after blowing up $B_0$,
we get one
extra base curve $C_j=E\cap H_j$
such that  the central fiber has multiplicity $3$ at $C_j$ and the fiber over the point $(1:0)$ has multiplicity $1$ at $C_j$.
The blow up of $E\cap H_0$ does not produce any new base curves. Thus, we have
$$
\frac{3}{1}-1=2
$$
exceptional divisors in the central fiber
lying over $P$, and in total~$3\cdot 2=6$ divisors in the central fiber lying over $P_{0j}$.

Now set $P=P_{ij}$ for some $1\leq i<j\leq 3$.
Resolving the base set of the pencil in a neighborhood of general point of
$B_i$, we have five exceptional di\-vi\-sors~$E_1,\ldots,E_5$ in the central fiber lying over $P$,
and the multiplicity of the central fiber at $E_s$  is $6-s$,
so that we have five new irreducible components
$B^s_j=E_s\cap H_j$.
The multiplicity of the central fiber at $B^s_j$ is $6-s$, while the multiplicity of the fiber over $(1:0)$ at $B^s_j$ is $1$.
Therefore, in total
we have
$$
\left(\frac{5}{1}-1\right)+\left(\frac{4}{1}-1\right)+\left(\frac{3}{1}-1\right)+\left(\frac{2}{1}-1\right)+\left(\frac{1}{1}-1\right)=10
$$
exceptional divisors in the central fiber lying over $P$,
and $3\cdot 10=30$ divisors lying over~$P_{ij}$.

Summing up, we have $1+16+6+30=53$ components of the central fiber of our log Calabi--Yau compactification for $f_X$, so Conjecture~\ref{conjecture:mirror Hodge}
holds for it, since $h^{12}(X)=52$.
Conjecture~\ref{conjecture:mirror Hodge} for $X$ is proven in~\cite{Prz13}, but the calculations we did to check it in fact are similar to ones suggested in~\cite{CP18}.
\end{example}

The general case can be done similarly, 
it just differs by more noisy combinatorics. That is, we need to count exceptional divisors lying over the strata
of base set of the initial family which are intersections of the multiple hyperplane with some number of hyperplanes with some multiplicities.
This approach was developed in~\cite{PS15}. Say, in the neighborhood of~$B_{ij}=B_i\cap B_j$ for $1\leq i<j\leq \alpha$
(we use the notation of the proof of Theorem~\ref{theorem:coverings}) the pencil looks like~$\{\mu z^{ad}=\lambda xy\}$, where $x,y,z$ are coordinates in $\Aff^{\alpha+1}$.
By~\cite[Proposition 4.7]{PS15}, the number of components of the central fiber lying over $B_{ij}$ is
$$\binom{ad-1}{2}=\frac{(ad-1)(ad-2)}{2}.$$

\begin{problem}
\label{problem:Hodge for coverings}
Prove Conjecture~\ref{conjecture:mirror Hodge} for toric Landau--Ginzburg models of Givental's type for Fano index one coverings of projective spaces.
The most promising approach is to generalize one suggested in~\cite{PS15}.
\end{problem}

\section{The general case}
\label{section:discussion}
Compactifications constructed in Sections~\ref{section:del Pezzo} and~\ref{section:coverings} suggest the following singular log Calabi--Yau compactification procedure
for weighted complete intersections. Consider a smooth Fano weighted complete intersection $X\subset \PP$ having a nice nef-partition. Let us use the notation of Definitions~\ref{definition:nef-partition} and~\ref{definition:Givental LG}. In particular, let us have a weak Landau--Ginzburg model of Givental's type $f_X$
corresponding to $X$.
Put $i_X=d_0=\sum_{j=0}^N a_j-\sum_{i=1}^k d_i$.
One can directly prove the following.

\begin{lemma}
The polytope $\nabla$ dual to the fan polytope $F(f_X)$ is generated by the rows of the matrix
$$
M=\left(%
\begin{array}{rrrr|r|rrrr|rrr}
  \frac{i_X}{a_{11}} & 0 & \ldots & 0 & \ldots & 0 & 0 & \ldots & 0 & -1 & \ldots & -1 \\
  0 & \frac{i_X}{a_{12}} & \ldots & 0 & \ldots & 0 & 0 & \ldots & 0 & -1 & \ldots & -1 \\
  \ldots & \ldots & \ldots & \ldots & \ldots & \ldots & \ldots & \ldots & \ldots & \ldots & \ldots & \ldots\\
  0 & 0 & \ldots & \frac{i_X}{a_{1m_1}} & \ldots & 0 & 0 & \ldots & 0 & -1 & \ldots & -1 \\
  -\frac{i_X}{a_{10}} & -\frac{i_X}{a_{10}} & \ldots & -\frac{i_X}{a_{10}} & \ldots & 0 & 0 & \ldots & 0 & -1 & \ldots & -1 \\
  \hline
  \ldots & \ldots & \ldots & \ldots & \ldots & \ldots & \ldots & \ldots & \ldots & \ldots & \ldots & \ldots\\
  \hline
  0 & 0 & \ldots & 0 & \ldots & \frac{i_X}{a_{k1}} & 0 & \ldots & 0 & -1 & \ldots & -1 \\
  0 & 0 & \ldots & 0 & \ldots & 0 & \frac{i_X}{a_{k2}} & \ldots & 0 & -1 & \ldots & -1 \\
  \ldots & \ldots & \ldots & \ldots & \ldots & \ldots & \ldots & \ldots & \ldots & \ldots & \ldots & \ldots\\
  0 & 0 & \ldots & 0 & \ldots & 0 & 0 & \ldots & \frac{i_X}{a_{km_k}} & -1 & \ldots & -1 \\
  0 & 0 & \ldots & 0 & \ldots & -\frac{i_X}{a_{k0}} & -\frac{i_X}{a_{k0}} & \ldots & -\frac{i_X}{a_{k0}}& -1 & \ldots & -1 \\
\hline
  0 & 0 & \ldots & 0 & \ldots & 0 & 0 & \ldots & 0 & \frac{i_X}{a_{0,1}}-1 & \ldots & -1 \\
  \ldots & \ldots & \ldots & \ldots & \ldots & \ldots & \ldots & \ldots & \ldots & \ldots & \ldots & \ldots \\
  0 & 0 & \ldots & 0 & \ldots & 0 & 0 & \ldots & 0 & -1 & \ldots & \frac{i_X}{a_{0,m_0}}-1 \\
\end{array}%
\right).
$$
\end{lemma}

\begin{remark}
Let us call a nef-partition $(I_1,\ldots,I_k)$ \emph{strong} if $a_{0,j}=1$ for $j=1,\ldots,m_0$
and $a_{i,j}$ divides $d_i$ for all $i,j$.
We expect that for any smooth Fano weighted complete intersection we always can find a strong nef-partition.
For instance, this holds for the cases listed in Theorem~\ref{theorem:nef-partitions}.
(Of course, a weighted complete intersection can have not strong nef-partition as well; the example is
a complete intersection of quartic and sextic in $\PP(1,1,1,1, 2,2,3)$.)
\end{remark}

Let $l$ be the number of weights of $\PP$ which are equal to $1$.
Assume for simplicity that~$i_X=1$, so that $m_0=0$ and we do not have the SE block of the matrix $M$.
By~\cite[Theorem 3.3.4]{Do82} and~\cite[Section 3.4.3]{Do82}, one has
$$
h^0(\mathcal O_X (-K_X))=h^0(\mathcal O_X(1))=l.
$$
The integral points of the boundary of $\nabla$ correspond to those rows of $M$ that have $a_{ij}=1$,
so there are $l-1$ of them. 
Let $T$ be the toric variety such that integral generators of rays of its fan have coordinates given by rows of~$M$.
Then
$
T=\PP^{m_1}\times\ldots\times \PP^{m_k}.
$
Compactify the family of fibers for $f_X$ in $T$ to a family $\mathcal F$.
Then $\mathcal F_\infty$ is a union of boundary divisors for $T$,
and their multiplicities are equal to $a_{ij}$.
The fiber $\mathcal F_0$ is a union of components of type
$$
\PP^{m_1}\times\ldots\times \PP^{m_{i-1}}\times H_i\times \PP^{m_{i+1}}\times P\times \PP^{m_k},
$$
where $H_i$ is a hyperplane section of $\PP^{m_i}$; multiplicity of such component is equal to $d_i$.
Let $D_{ij}$ be the boundary divisor that correspond to the row of $M$ associated with $a_{ij}$.
Note that
$$
\mathcal F_\lambda\sim -K_{T}+\sum (a_{ij}-1)D_{ij}.
$$

So we suggest the following compactification procedure.
Choose a strong nef par\-ti\-ti\-on~$\overline I=(I_1,\ldots,I_k)$  and a weak Landau--Ginzburg model $f_X$ associated with it.
Compactify~$f_X$ in $T$ and get the pencil $\mathcal F$.
Let $D_1,\ldots,D_r$ be the boundary divisors of $T$ whose multiplicities $s_i$ in the base locus is greater then $1$, and let~$B_i\subset D_i$ be the component of the base set lying on $D_i$.
  Resolve the base set in a neighbourhood of $B_i$, $i=1,\ldots,r$, by $s_i$ blow ups.
  Note that since $\overline I$ is strong, all exceptional divisors, except for the last horizontal one,
  lie in the central fiber. 
Now make toric birational transformation in codimension $2$ if needed and blow down proper transforms of $D_i$.
We get a family $\mathcal F'=\{\mathcal F'_t,\ t\in \PP^1\}$ with toric total space $Z'$ such that  $\mathcal F'_t\sim -K_{Z'}$. Note that the base set of $\mathcal F'$  is a union of smooth components of codimension $2$, and these components do not coincide with
  intersections of components of the fiber over infinity.
  The variety $Z'$ can be singular; the singularities should not intersect the base set of $\mathcal F'$.
  Now we can proceed as in Compactification Construction~\ref{compactification construction}, cf.~\cite[Proposition 26]{Prz17},
and get the required singular log Calabi--Yau compactification $Z\to \PP^1$.
Note that if the suggested compactification procedure gives a (singular) log Calabi--Yau compactification,
then Conjecture~\ref{conjecture:linear systems} holds for $X$.
Note also that if $Z$ is singular, then in the same way as in the proof of Lemma~\ref{lemma: non-projective} one can prove that there is no smooth projective
log Calabi--Yau compactification of $f_X$.

The case $i_X>1$ can be considered in the similar way, cf.~\cite{Prz18} and ~\cite[\S 9.3]{PSh}.

\begin{remark}
If we choose a non-strong nef-partition, we get exceptional divisors at the fiber over infinity, so we need to prove that we can contract them.
\end{remark}

\begin{problem}
Prove that the suggested compactification procedure works. 
\end{problem}

\end{document}